\documentclass{amsart}

\usepackage{amsfonts}
\usepackage{amsmath}
\usepackage{amsthm}
\usepackage{amssymb}
\usepackage[english]{babel}
\usepackage{graphics}
\usepackage{latexsym}
\usepackage{longtable} \usepackage{mathrsfs}
\usepackage{graphicx}
\usepackage{wrapfig} 
\usepackage[pdftex,colorlinks=true]{hyperref}

\newtheorem{theorem}{Theorem}
\newtheorem{corollary}[theorem]{Corollary}
\newtheorem{lemma}[theorem]{Lemma}
\newtheorem{proposition}[theorem]{Proposition}

\newtheorem{claim}[theorem]{Claim}
\newtheorem{example}[theorem]{Example}

\theoremstyle{definition}
\newtheorem{definition}[theorem]{Definition}
\newtheorem{remark}[theorem]{Remark}

\renewcommand{\S}{\mathcal{S}}

\newcommand{\R}{\mathbb{R}}
\newcommand{\N}{\mathbb{N}}
\newcommand{\mS}{\mathbb{S}}
\newcommand{\mB}{\mathbb{B}}

\newcommand{\noi}{\noindent}
\newcommand{\ms}{\medskip}

\newcommand{\al}{\alpha}
\newcommand{\be}{\beta}
\newcommand{\ga}{\gamma}
\newcommand{\de}{\delta}
\newcommand{\De}{\Delta}
\newcommand{\e}{\varepsilon}
\newcommand{\si}{\sigma}
\newcommand{\la}{\lambda}

\newcommand{\Om}{\Omega}

\newcommand{\larrow}{\longrightarrow}
\newcommand{\ot}{\otimes}

\newcommand{\ri}{\rightarrow}

\newcommand{\p}{\partial}
\newcommand{\sub}{\subseteq}
\newcommand{\set}{\setminus}
\newcommand{\by}{\times}
\newcommand{\rk}{\textrm{rk}}

\newcommand{\tr}{\textrm{tr}}
\newcommand{\sgn}{\textrm{sgn}}
\newcommand{\ess}{\textrm{ess}}

\newcommand{\dist}{\textrm{dist}}

\newcommand{\Div}{\textrm{Div}}

\newcommand{\inter}{\textrm{int}}

\newcommand{\cof}{\textrm{cof}}

\newcommand{\bt}{\begin{theorem}}\newcommand{\et}{\end{theorem}}
\newcommand{\bd}{\begin{definition}}\newcommand{\ed}{\end{definition}}
\newcommand{\bl}{\begin{lemma}}\newcommand{\el}{\end{lemma}}
\newcommand{\beq}{\begin{equation}}\newcommand{\eeq}{\end{equation}}
\newcommand{\bc}{\begin{claim}}\newcommand{\ec}{\end{claim}}
\newcommand{\bex}{\begin{example}}\newcommand{\eex}{\end{example}}
\newcommand{\bcor}{\begin{corollary}}\newcommand{\ecor}{\end{corollary}}
\newcommand{\bp}{\begin{proof}}\newcommand{\ep}{\end{proof}}

\newcommand{\BPL}{\medskip \noindent \textbf{Proof of Lemma} }

\newcommand{\BPCOR}{\medskip \noindent \textbf{Proof of Corollary} }
\newcommand{\BPP}{\medskip \noindent \textbf{Proof of Proposition} }

\numberwithin{equation}{section}
\numberwithin{theorem}{section}

\begin{document}

\title{Optimal $\infty$-Quasiconformal Immersions}

\author{\textsl{Nikos Katzourakis}}
\address{Department of Mathematics and Statistics, University of Reading, Whiteknights, PO Box 220, RG6 6AX, UK and BCAM, Alameda de Mazarredo 14, E-48009, Bilbao, Spain}
\email{n.katzourakis@reading.ac.uk}


\date{}



\begin{abstract} For a Hamiltonian $K \in C^2(\R^{N \by n})$ and a map  $u:\Om \sub \R^n \!\larrow \R^N$, we consider the supremal functional
\[  \label{1} \tag{1}
E_\infty (u,\Om) \ :=\ \big\|K(Du)\big\|_{L^\infty(\Om)}  .
\]
The ``Euler-Lagrange" PDE associated to \eqref{1} is the quasilinear system
 \[  \label{2}
A_\infty u \, :=\, \Big(K_P \ot K_P +  K[K_P]^\bot \! K_{PP}\Big)(Du):D^2 u \, =  \, 0. \tag{2}
\]
Here $K_P$ is the derivative and $[K_P]^\bot$ is the projection on its nullspace. \eqref{1} and \eqref{2} are the fundamental objects of vector-valued Calculus of Variations in $L^\infty$ and first arose in recent work of the author \cite{K1}-\cite{K6}. Herein we apply our results to Geometric Analysis by choosing as $K$  the dilation function 
\[
K(P)={|P|^2}{\det(P^\top\! P)^{-1/n}} 
\]
which measures the deviation of $u$ from being conformal. Our main result is that appropriately defined minimisers of \eqref{1} solve \eqref{2}. Hence, PDE methods can be used to study optimised quasiconformal maps. Nonconvexity of $K$ and appearance of interfaces where $[K_P]^\bot$ is discontinuous cause extra difficulties. When $n=N$, this approach has previously been followed by Capogna-Raich \cite{CR} and relates to Teichm\"uller's theory. In particular, we disprove a conjecture appearing therein. 
\end{abstract}

\maketitle

\section{Introduction} \label{section1}

Let $\mathcal{M}_0$ be a topological submanifold of $\R^N$ with boundary. In this paper we are interested in the problem of finding a Riemannian manifold $(\mathcal{M},g)$ which has \emph{minimal dilation} and satisfies $\p \mathcal{M} =\p \mathcal{M}_0$. In this setting, dilation is a functional on $L^\infty(\mathcal{M},\ot^{(2)} T^*\mathcal{M})$, defined as the $L^\infty$ norm of the trace of the Distortion Tensor
\beq \label{1.1}
\bold{G}\ :=\ \frac{g}{\det(g)^{1/n}} .
\eeq
This problem is an extension of the classical \emph{Teichm\"uller Problem} (see \cite{T, AIM, AIMO}). The scaling in \eqref{1.1} is such that $\bold{G}$ is invariant under conformal tranformations and, as we explain later, the geometric meaning of $\tr(\bold{G})$ being ``minimal" is that ``geometry is distorted as less as possible''. As a first step, we consider a simplified problem for the case of immersions $u : \Om \sub \R^n \larrow \R^N$ with prescribed boundary values on $\p \Om$. Then,  the dilation functional for $\mathcal{M}=u(\Om)$ becomes
\beq \label{1.2}
K_\infty(u,\Om)\ :=\ \big\|K(Du)\big\|_{L^\infty(\Om)},
\eeq 
where $K$ will be called the \emph{dilation function} and is given by
\beq \label{1.3}
K(P)\ :=\ \left\{\begin{array}{l}
           \dfrac{|P|^2}{\det(P^\top \! P )^{1/n}}\ ,  \ \ \ \text{ on }S^+,\ms\\
             +\infty \ ,   \ \ \ \ \ \ \ \ \ \ \ \ \ \ \ \, \text{ on }\R^{N \by n} \set S^+.
           \end{array}
\right.
\eeq
In \eqref{1.3}, $|P|=\tr(P^\top\! P)^{1/2}$ is the Euclidean norm on $\R^{N \by n}$ and 
\beq
S^+ \ :=\ \Big\{P \in \R^{N \by n} \ :\  \det\big(P^\top \! P \big)>0\Big\}. 
\eeq
Important objects of Geometric Topology related to \eqref{1.2} arise for $n=N$. Homeomorphisms $u : \Om \sub \R^n \larrow \R^n$ in $W^{1,n}_{loc}(\Om)^N$ which satisfy $K_\infty(u,\Om)<\infty$ are called \emph{Quasiconformal Maps} and constitute a class of maps well studied in the literature; see for example \cite{Ah2, B, G, S, V}. $L^p$ averages of Quasiconformal maps, that is weakly differentable homeomorphisms for which  $\| K(Du) \|_{L^p(\Om)}<\infty$ have also been systematically considered. \emph{Conformal maps}, namely those homeomorphisms $u : \Om \sub \R^n \larrow \R^n$  in $C^1(\Om)^N$ which satisfy $Du^\top Du = \frac{1}{n}|Du|^2I$ form a special important class of Quasiconformal maps since for those $K(Du)$ is constant and equals $n$. Conformal maps preserve \emph{angles}, but not necessarily \emph{lengths} and hence distort the geometry of $\Om$ in a controlled fashion.  However, by Liouville's rigidity theorem, when $n\geq 3$ the only conformal maps that exist are compositions of rotations, dilations, and the inversion $x \mapsto x/|x|^2$. Hence, quasiconformal maps for which $K(Du)$ is merely bounded relax conformality but still deform $\Om$ to $u(\Om)$ in a fairly controlled fashion.

The problem with Quasiconformal maps is that too little information on their structure is provided by a mere norm bound, and the same holds for the \emph{finite distortion problem} when one restricts attention to minimisers of the dilation functional. The subtle point is that \eqref{1.2} is \emph{nonlocal}, in the sense that with respect to the $\Om$ argument \eqref{1.2} is not a measure. Simple examples certify that minimisers over a domain with fixed boundary values are not local minimisers over subdomains and the direct method of Calculus of Variations when applied to \eqref{1.2} generally does not produce PDE solutions. 

In the very recent work, Capogna and Raich \cite{CR}, remedied this problem by ``optimising" Quasiconformal maps. The idea is to consider an appropriate nonstandard $L^\infty$ variational problem for \eqref{1.2} and derive a PDE governing Optimal Quasiconformal Maps that can be used as platform for their qualitative study. Motivated by the classical results of Aronsson \cite{A1, A2} on \emph{Calculus of Variations in $L^\infty$}, they developed an $L^\infty$ variational approach for extremal (as they are called therein) quasiconformal maps. The essence of this approach is the following: let $Q_p u =0$ be the Euler-Lagrange system associated to the functional $\|K(Du)\|_{L^p(\Om)} $. Then, at least formally $Q_p$ tends to a certain operator $Q_\infty$ and  $\|K(Du)\|_{L^p(\Om)}$ tends to $\|K(Du)\|_{L^\infty(\Om)}$, both as $p\ri \infty$. The operator $Q_\infty$ defines a quasilinear 2nd order system in non-divergence form. However, it is not a priori clear that the following rectagle ``commutes" 
\begin{align} \label{1.8}
&\|K(Du)\|_{L^p(\Om)}   \ \ \ \  \larrow\ \ \ \ Q_p u =0  \nonumber\\
& \ \ \ \downarrow\ p \ri \infty\ \ \ \ \ \ \ \ \ \ \  \ \ \ \ \ \ \downarrow\ p \ri \infty  \\
&\|K(Du)\|_{L^\infty(\Om)}   \ \ \ \ \dashrightarrow \ \ \ \ Q_\infty u=0   \nonumber
\end{align}
so that $Q_\infty$ has a variational structure with respect to $K_\infty$, in the sense that appropriately defined minimisers of $K_\infty$ $u: \Om \sub \R^n \larrow \R^n$ solve $Q_\infty u =0$.  In such an event, $Q_\infty u=0$ will play the role of ``Euler-Lagrange PDE" for the dilation functional. This turns out to be the case, though. Among other far-reaching contributions which include a deep study of dilations of  extensions up to the boundary and quasiconformal gradient flows, Capogna and Raich introduced in \cite{CR} a localized minimality notion for \eqref{1.2} and proved that those local minimisers  among ``competitors" indeed solve the formally derived PDE.

Simultaneously and independently, the author, also inspired by Aronsson's work and the successful modern evolution of the field of Calculus of Variations in $L^\infty$ (see for example \cite{C}), initiated the development of vector-valued Calculus of Variations in $L^\infty$ for general supremal functionals in \cite{K1}-\cite{K6} with particular emphasis to the model functional $\|Du\|_{L^\infty(\Om)}=\ess \, \sup_\Om |Du|$. For a Hamiltonian $H\in C^2(\R^{N \by n})$ and the respective supremal functional
\beq \label{1.6a}
E_\infty(u,\Om)\ :=\ \|H(Du)\|_{L^\infty(\Om)},
\eeq
the PDE system which plays the role of ``Euler-Langrange PDE" for \eqref{1.6a} is
\beq  \label{1.4a}
A_\infty u \, :=\, \Big(H_P \ot H_P + H[H_P]^\bot  H_{PP}\Big)(Du):D^2 u \, =  \, 0 .
\eeq
Here $[H_P(Du(x))]^\bot$ is the projection on the nullspace of $H_P(Du(x))^\top : \R^N \larrow \R^n$,  
and $H_P,H_{PP}$ denotes derivatives (for details see Preliminaries \ref{section2}). The special  case of $H(P)=|P|^2$ leads to the important \emph{$\infty$-Laplacian}
\beq  \label{1.5a}
\De_\infty u \, :=\, \Big(Du \ot Du +  |Du|^2[Du]^\bot \! \ot I\Big):D^2 u \, =  \, 0. 
\eeq
System \eqref{1.4a} is a quasilinear 2nd order system in non-divergence form which arises in the limit of the Euler-Lagrange system of the $L^p$ functional $\|H(Du)\|_{L^p(\Om)}$ as $p\ri \infty$. In the scalar case of $n=1$ the normal coefficient of \eqref{1.5a} $ |Du|^2[Du]^\bot$ vanishes, and the same holds for submersions in general. The scalar $\infty$-Laplacian then becomes $Du\ot Du :D^2u=0$.

Unlike the scalar case of $n=1$, in the full vector case of \eqref{1.4a} intriguing phenomena appear. Except for the emergence of ``singular solutions'' to \eqref{1.4a}, a further difficulty not present in the scalar case is that \emph{\eqref{1.4a} has discontinuous coefficients} even for $C^\infty$ solutions. There exist $C^\infty$ solutions whose rank of $H_P(Du)$ is not constant: such an example on $\R^2$ for \eqref{1.5a} is given by $u(x,y) = e^{ix}-e^{iy}$ which is $\infty$-Harmonic near the origin and has $\rk(Du)=1$ on the diagonal, but it has $\rk(Du)=2$ otherwise and hence the projection $[Du]^\bot$ is discontinuous (\cite{K1}). More sophisticated examples with interfaces which have junction and corners appear in \cite{K4}. In general, \emph{$\infty$-Harmonic maps present a phase separation} and on each phase the dimension of the tangent space is constant and these phases are separated by \emph{interfaces} whereon the rank of $Du$ ``jumps'' and $[Du]^\bot$ is discontinuous (\cite{K1}, \cite{K6}). Extensions of the results of \cite{K1}, \cite{K2} to the subelliptic setting appear in \cite{K3}. Moreover, it has very recently been established that the celebrated scalar $L^\infty$ uniqueness theory has no counterpart when $N\geq 2$ (\cite{K5}).

In this paper we work towards the problem mentioned in the beginning by extending the theory of \cite{CR} to the case of immersions $u : \Om \sub \R^n \larrow \R^N$ and in the same time we elaborate it and make it more efficient in certain respects. First of all, we allow for positive codimension $N-n$ and take into account the exterior geometry of immersions. Moreover, our maps are local diffeomorphisms onto their images, but in our analysis we do \emph{not} impose the global topological constraint that our maps are homemorphisms onto their image and allow for self-intersections. However, \emph{all} our results and notions are still valid and with the exact same proofs in this restricted class. For distinction, we introduce the following terminology: \emph{an immersion $u : \Om \sub \R^n \larrow \R^N$ in $C^1(\Om)^N$ is called $p$-Quasiconformal when $\| K(Du) \|_{L^p(\Om)}<\infty$, $1\leq p \leq \infty$}. We begin by repeating part of the program of \cite{K1}, \cite{K2} under the lens of \cite{CR} to the extended case. After some introductory material is Section \ref{section2}, in Section \ref{section3} we calculate the PDE system which Optimal $p$-Quasiconformal immersions satisfy (equations \eqref{3.8}, \eqref{3.9}), that is the Euler-Lagrange system of $K_p(u,\Om):= \| K(Du) \|_{L^p(\Om)}$. Then, in Section \ref{section4} we formally derive in the limit as $p\ri \infty$ the PDE system which Optimal $\infty$-Quasiconformal immersions $u: \Om \sub \R^n \larrow \R^N$ satisfy, that is the system associated to \eqref{1.2}:
\beq  \label{1.4}
Q_\infty u \, :=\, \Big(K_P \ot K_P + K[K_P]^\bot K_{PP}\Big)(Du):D^2 u \, =  \, 0
\eeq
where the derivatives of the dilation are given by
\begin{align} 
K_P(Du)\ &=\ 2Du \frac{g^{-1}S(g)}{\det(g)^{1/n}} \label{2.10a},\\
K_{PP}(Du)\ &= \ 2\left(I \ot  \frac{g^{-1}S(g)}{\det(g)^{1/n}} \ + \, Du \ot Du : \frac{g^{-1}E}{\det(g)^{1/n}}\right)\ + \ O(Du). \label{2.10b}
\end{align}
Here $g=Du^\top\! Du$, $S$ is the Ahlfors operator given by \eqref{2.6}, $E$ is a constant tensor given by \eqref{3.3A} and $O(Du)$ is a tensor annihilated by $[K_P(Du)]^\bot$ and does not appear in the PDE system \eqref{1.4} (for details see Lemmas \ref{l1},  \ref{l1,1}). The derivation has overlaps with the respective in \cite{K1}, but is not a direct consequence since we utilise the specific structure of the Hamiltonian \eqref{1.3}. By restricting ourselves to $n=N$ and employing Lemma \ref{l6} to relate the seemingly different system \eqref{1.4} to that of \cite{CR}, we see that the derivation as $p\ri \infty$ in \cite{CR} is incomplete and their PDE is only a part of \eqref{1.4}. System \eqref{1.4} consists of two systems whose defining vector-valued nonlinearities are normal to each other:
\begin{align}
K_P(Du) \ot K_P(Du) :D^2u\, =&\, 0, \label{1.11}\\
[K_P(Du)]^\bot K_{PP}(Du): D^2 u \, =&\, 0.  \label{1.12}
\end{align}
System \eqref{1.11} is the ``tangential" part in (the range of the projection) $[K_P(Du)]^\top$ and system \eqref{1.12} is the ``normal" part in $[K_P(Du)]^\bot$ (see Figure 1). The reason for this terminology is that $[Du]^\top$ is (the projection on) the tangent bundle of the immersion, $[Du]^\bot$ is its normal bundle and by \eqref{2.10a} we have that $[K_P(Du)]^\top \sub [Du]^\top$.

\[
\underset{\text{Figure 1.}}{\includegraphics[scale=0.24]{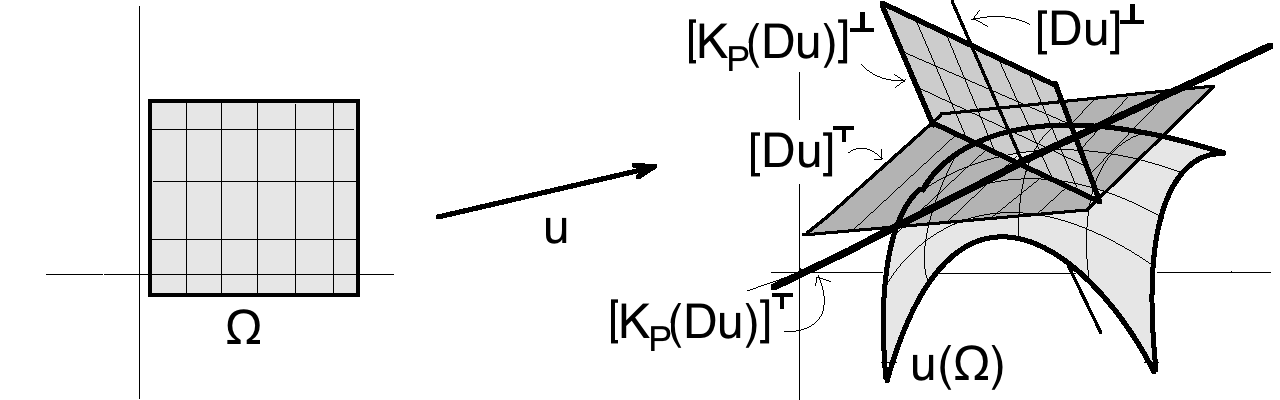}} \label{fig1} 
\]
The derivation in \cite{CR} has lost information along directions in $[K_P(Du)]^\bot$ and reveals only system \eqref{1.11}. System \eqref{1.12} appears also in zero-codimension when $n=N$ since generally $K_P(Du)$ does not have rank equal to $n$,  although by assumption the rank of $Du$ equals $n$. More importantly, \emph{when the rank of $K_P(Du)$  becomes nonconstant, the coefficients of \eqref{1.4} become discontinuous}. This leads to the \emph{appearance of interfaces} whereon the projection $[K_P(Du)]^\bot$  is discontinuous. These interfaces are boundaries of  the different \emph{phases} to which Optimal $\infty$-Quasiconformal maps naturally separate.

In Section \ref{section5} we move to the variational structure of Optimal $\infty$-Quasiconformal maps. Inspired from \cite{K2}, we introduce the variational notion of \emph{$\infty$-Minimal Dilation}, which is Rank-One Locally Minimal Dilation with ``Minimally Distorted Area'' of $u(\Om)$ (Definition \ref{def1}). Rank-one locally minimal dilation requires that an immersion is a local minimiser for the dilation functional when the ``set of competitors" is the one obtained by taking essentially scalar local variations with fixed zero boundary values (Figure 2). Minimally distorted area means that the immersion is a local minimiser where the ``set of competitors" is the one obtained by taking variations along sections of the normal  vector bundle $[K_P(Du)]^\bot$ over $u(\Om)$ with free boundary values (Figure 3). The appearance of interfaces where the dimension of $[K_P(Du)]^\bot$ jumps causes substantial difficulties, even in the very definition of the minimality notion. Our first main result is Theorem \ref{th1}, wherein we prove that $\infty$-Quasiconformal maps with $\infty$-Minimal Dilation are Optimal, \emph{at least off the interfaces of discontinuities in the coefficients}. This result follows closely Theorem 2.1 in \cite{K1} and Theorem 2.2 in \cite{K2}, but nonconvexity of \eqref{1.3}, appearance of discontinuities in \eqref{1.4} and the necessity of restriction to specific variations create complications not present in the results just quoted. We note that the rank-one minimality notion gives rise to the tangential system and the condition on the minimality of the area gives rise to the normal system. 

In Section \ref{section6} we study some geometric aspects of \eqref{1.4} and of the interfaces of its solutions. In Subsection \ref{subsection6.1} we show that system \eqref{1.4} has a ``geometric" rather coordinate-free reformulation, at least off interfaces of discontinuities. More precisely, \eqref{1.11} and \eqref{1.12} are respectively equivalent to
\begin{align} 
S(\bold{G})D\big(\tr(\bold{G})\big)\ &= \ 0 , \label{1.13}\\ 
\mathbb{B}^\bot : \big(\tr(\bold{G})\big)_{P} \ &=  \ 0,\label{1.14}
\end{align}
where $\bold{G}$ is given by \eqref{1.1} for $g=Du^\top \! Du$ and $\mathbb{B}^\bot$ is a ``generalized 2nd fundamental form" with respect to normal sections valued in $[K_P(Du)]^\bot$. If $K_P(Du)$ has full rank $n$, then $[K_P(Du)]^\bot$ coincides with the normal bundle $[Du]^\bot$ of the immersion and $\mB^\bot$ reduces to the standard object. System \eqref{1.13} is quite ``metrically invariant" but system \eqref{1.14} depends on the exterior geometry and measures the ``shape of $u(\Om)$". In Subsection \ref{subsection6.2}, by assuming some a priori local $C^1$ regularity on the interfaces but with possible self-intersections, we prove an identity which shows that the covariant gradient of $[K_P(Du)]^\bot$ along the interface is differentiable when projected along $K_P(Du)$. 

In Section \ref{section7} we turn our attention to the converse statement of that in Theorem \ref{th1}, that is the sufficiency of \eqref{1.4} for the variational notion of $\infty$-Minimal Dilation. Nonconvexity of \eqref{1.3} and the resemblance to similar phenomena in \emph{Minimal Surfaces} leaves little hope for system  \eqref{1.12} to be sufficient for minimally distorted area. However, in Proposition \ref{c8} we establish that when $n=2\leq N$ there is a triple equivalence among solutions of \eqref{1.11}, the condition the dilation \eqref{1.3} to be constant and the immersion to have rank-one locally minimal dilation. This result relates directly to the two-dimensional results in \cite{Ah1, B, H}. In particular, when $n=2$ interfaces disappear and the coefficients of \eqref{1.4} become continuous. 

Moreover, as a consequence of Example \ref{ex1} which certifies that rank-one locally minimal dilation is \emph{strictly weaker} than the variational notion utilized in \cite{CR} with respect to general vector-valued variations (among competitors), \emph{we disprove the conjecture of Capogna-Raich on the sufficiency of \eqref{1.3} explicitely stated in p.\ 855}. Finally, at the end of Section \ref{section7} we loosely discuss the much more complicated case when $n\geq 3$. In this case results are less sharp. Although it is hardly conclusive, it seems that dilation may not be constant but we do believe that \eqref{1.11} is still sufficient for rank-one locally minimal dilation.

Throughout this paper, as in \cite{CR} and also in \cite{K1}-\cite{K6}, we restrict our analysis to the unnatural class of $C^2$ solutions. This is only the first step in our study and we can not go much further without an appropriate ``weak'' theory of nondifferentiable solutions for \eqref{1.4}. The latter  much deeper problem, namely defining a notion of solution for which we can also prove existence to the Dirichlet problem, will be considered in future work.

\section{Preliminaries.} \label{section2} 

Throughout this paper we reserve $n,N \in \N$ for the dimensions of Euclidean spaces and $\mS^{N-1}$ denotes the unit sphere of $\R^N$. Greek indices $\al, \be, \ga,... $ run from $1$ to $N$ and Latin $i,j,k,...$ form $1$ to $n$. The summation convention will always be employed in repeated indices in a product. Vectors are always viewed as columns and we differentiate along rows. Hence, for $a,b\in \R^n$, $a^\top b$ is their inner product and $ab^\top$ equals $a \ot b$. If $u=u_\al e_\al :  \Om \sub \R^n \larrow \R^N$ is in $C^2(\Om)^N$, the gradient matrix $Du$ is viewed as $D_i u_\al e_\al \ot e_i : \Om \larrow \R^{N \by n}$ and the Hessian tensor $D^2u$ as $D^2_{ij} u_\al e_\al \ot e_i \ot e_j: \Om \larrow \R^{N \by n^2}$. The Euclidean (Frobenious) norm on $\R^{N\by n}$ is $|P|=(P_{\al i}P_{\al i})^{1/2} = (\tr (P^\top P))^{1/2}$.  We also introduce the following \emph{contraction operation} for tensors which extends the Euclidean inner product $P:Q=\tr(P^\top Q)=P_{\al i}Q_{\al i}$ of $\R^{N\by n}=\R^N \ot \R^n$. Let ``$\ot^{(r)}$'' denote the $r$-fold tensor product. If $S\in \ot^{(q)}\R^N \ot^{(s)} \R^n$,  $T \in \ot^{(p)}\R^N \ot^{(s)} \R^n$ and $q\geq p$, we define a tensor $S:T$ in $\ot^{(q-p)} \R^N$ by 
\beq \label{1.24}
S:T \ :=\ \big(S_{\al_q ...\al_p... \al_1 \, i_s ... i_1}  T_{\al_{p}  ... \al_1 \, i_s ... i_1}  \big) \, e_{\al_q} \ot ... \ot e_{\al_{p+1}}. 
\eeq
For example, for $s=q=2$ and $p=1$, the tensor $S:T$ of \eqref{1.24} is a vector with components $S_{\al \be i j}T_{\be ij}$ with free index $\al$ and  the indices $\be,i,j$ are contracted. In particular, in view of \eqref{1.24}, the 2nd order linear system 
\beq
A_{\al i \be j}D^2_{ij}u_\be \, +\, B_{\al \ga k} D_ku_\ga + C_{\al \de} u_\de\, =\, f_\al ,
\eeq
can be compactly written as $A$:$D^2u + B$:$Du+Cu=f$, where the meaning of ``$:$" in the respective dimensions is made clear by the context. Let now $P : \R^n \larrow \R^N$ be linear map. The identity map of $\R^N$ splits as $I=[P]^\top \oplus [P]^\bot$, where $[P]^\top$ and $[P]^\bot$ denote orthogonal projection on range $R(P)$ and nullspace $N(P^\top)$ respectively. Moreover, for the dilation function \eqref{1.3}, we have $K(P)\geq n$ and $K(P)=n$ if and only if $P^\top\!P = \la I$ with $\la=\frac{1}{n}|P|^2$. This property of $K$ is a simple consequence of the inequality of arithmetic-geometric mean applied to the $n$ eigenvalues of $P^\top\!P$ by utilising the Spectral Theorem. Let us now recall some elementary properties of determinants. If $A=A_{ij}e_i \ot e_j \in \R^n \ot \R^n$, we have
\begin{align}
\cof(A)_{ij}\ :=& \ (-1)^{i+j} \det \Big(\underset{k\neq i, l \neq j}{\Sigma} A_{kl}e_k \ot e_l \Big),\\
\cof(A)\ :=& \ \cof(A)_{ij} e_i \ot e_j ,
\end{align}
\begin{align}
A\, \cof(A)^\top \ =& \ \cof(A)^\top A\ = \ \det(A)I,\\
D_{A_{ij}}\big(\det(A)\big)\, \equiv & \ \, \big(\det(A)\big)_{A_{ij}} = \ \cof(A)_{ij}.
\end{align}
Obviously, subscript denotes partial derivative. The \emph{Ahlfors operator} is defined by
\beq \label{2.6}
S(A)\ := \ \frac{1}{2}\big(A+A^\top \big)\ - \ \frac{1}{n}\tr(A) I
\eeq
and has the property that for any $A$, $S(A)$ is symmetric and  traceless, that is $\tr(S(A))=0$. Let now $u : \Om \sub \R^n \larrow \R^N$ be an immersion in $C^1(\Om)^N$. Then, the rank of $Du$ satisfies $\rk(Du)=n\leq N$. $u$ is \emph{Conformal} when there is $f\in C^0(\Om)$ such that $Du^\top \! Du =f^2 I$ on $\Om$, that is $D_i u_\al D_j u_\al=f^2\de_{ij}$.  For immersions, the Riemannian metric on $u(\Om)$ induced from $\R^N$ is  $g:= Du^\top \! Du$ and $g^{-1}$ denotes the pointwise inverse of the positive symmetric tensor $g$. Since $S(g)=g-\frac{1}{n}\tr(g)I$, we have the commutativity relation
\beq
g^{-1} S(g)\ =\ S(g)g^{-1}  \ =\  I\, -\,  \dfrac{\tr(g)}{n } g^{-1}
\eeq
which will be tacitly used in the sequel. In view of these conventions, the PDE system describing Optimal Quasiconformal immersions in index form reads
\beq  
\Big( K_{P_{\al i}} K_{P_{\be j}}  + \, K[K_P]_{\al \ga}^\bot K_{P_{\ga i}P_{\be j}} \Big)(Du) \, D^2_{ij} u_\be\ = \ 0.
\eeq
The derivatives $K_P,K_{PP}$ of $K$ appearing here and in \eqref{2.10a}, \eqref{2.10b} are given in index form by \eqref{3.1a}, \eqref{3.2A}. Finally, we will use the notation ``$\Gamma$'' for sections of vector bundles. We note that our terminology of ``$p$-Quasiconformal" slightly deviates from the usage of this term in the literature, but its purpose is to avoid the less elegant term ``$L^p$-Quasiconformal". Since we are only interested in the extreme case of $p=\infty$, there will be no confusion. We conclude by observing thatwhen $\Om \Subset \R^n$, all immersions $u : \overline{\Om} \sub \R^n \larrow \R^N$ in $C^1(\overline{\Om})^N$ are $p$-Quasiconformal for all $p\in [1,\infty]$.

\section{Derivation of the Euler-Lagrange PDE System Governing Optimal $p$-Quasiconformal Immersions.} \label{section3}

In this section we calculate the specific form of the Euler-Lagrange system associated to the functional $\| K(Du)\|^p_{L^p(\Om)}$ which Optimal $p$-Quasiconformal immersions satisfy. We begin by calculating first and second derivatives of \eqref{1.3}.

\begin{lemma} \label{l1} Let $K$ be given by \eqref{1.3}. Then, $K\in C^1(S^+)$ and its derivative is given by
\beq \label{3.1}
K_P(P)\ =\ 2P\frac{\big(P^\top \! P \big)^{-1}S\big(P^\top \! P \big)}{\det\big(P^\top \! P \big)^{1/n}}.
\eeq
\end{lemma}
In index form \eqref{3.1} can be written as
\beq \label{3.1a}
K_{P_{\al i}} (P)\ =\ 2P_{\al m} \Bigg( \frac{ \de_{mi}-\frac{1}{n}|P|^2\big(P^\top \! P \big)_{mi}^{-1} }{\det\big(P^\top \! P \big)^{1/n}} \Bigg).
\eeq

\BPL \ref{l1}. We begin by observing the triviality that for $P\in S^+$, the matrix $P^\top \! P$ is positive symmetric on $\R^n$ and also
\beq
\big(P^\top \! P \big)^{-1,\top}\ =\ \big(P^\top \! P \big)^{\top,-1}\ =\ \big(P^\top \! P \big)^{-1}.
\eeq
By differentiation of \eqref{1.3}, we have
\begin{align} \label{3.3}
K_{P_{\al i}}(P)\ =&  \ \frac{\ 2P_{\al i}\det\big(P^\top \! P \big)^{\frac{1}{n}} \, -\, \dfrac{|P|^2}{n} \det\big(P^\top \! P \big)^{\frac{1}{n}-1} \cof\big(P^\top \! P \big)_{kl}(P_{\be k}P_{\be l})_{P_{\al i}} }{\det\big(P^\top \! P \big)^{2/n}} 
\nonumber\\
    =&  \ \frac{\ 2P_{\al i} \, -\, \dfrac{|P|^2}{n \det\big(P^\top \! P \big) }\cof\big(P^\top \! P \big)_{kl} \big( \de_{\al \be} \de_{ik}P_{\be l} \, +\,  \de_{\al \be} \de_{il}P_{\be k} \big)}{\det\big(P^\top \! P \big)^{1/n}}.
\end{align}
Thus,
\begin{align} \label{3.4}
K_{P_{\al i}}(P)\ =&  \  \ \frac{\ 2P_{\al i} \, -\, \dfrac{|P|^2}{n \det\big(P^\top \! P \big) } \Big( \cof\big(P^\top \! P \big)_{il}  P_{\al l}  +\,  \cof\big(P^\top \! P \big)_{ki}  P_{\al k} \Big)}{\det\big(P^\top \! P \big)^{1/n}}
 \nonumber \\
        =&  \ 2P_{\al m} \, \frac{\ \de_{m i}\, -\, \dfrac{|P|^2}{n \det\big(P^\top \! P \big) } \dfrac{1}{2}\Big( {\cof\big(P^\top \! P \big)_{im}+ \, \cof\big(P^\top \! P \big)_{mi} } \Big)}{\det\big(P^\top \! P \big)^{1/n}}.
\end{align}
Hence, \eqref{3.4} gives
\beq \label{3.5}
K_P(P)\ = \ \frac{2P}{\det\big(P^\top \! P \big)^{1/n}}\left(I\, -\,  \dfrac{|P|^2}{n } \left( \dfrac{\cof\big(P^\top \! P \big)^\top +\,  \cof\big(P^\top \! P \big)}{2 \det\big(P^\top \! P \big)} \right)\right)
\eeq
and by using that 
\beq
\cof\big(P^\top \! P \big)^\top =\ \cof\big(P^\top \! P \big) \ =\ \big(P^\top \! P \big)^{-1}\det\big(P^\top \! P \big), 
\eeq
equation \eqref{3.5} gives
\begin{align} \label{3.6}
K_P(P)\ =& \ \frac{2P}{\det\big(P^\top \! P \big)^{1/n}}\left(I\, -\,  \dfrac{|P|^2}{n } \big(P^\top \! P \big)^{-1} \right)
\nonumber\\
     =& \ 2P \frac{\big(P^\top \! P \big)^{-1} }{\det\big(P^\top \! P \big)^{1/n}}\left(P^\top \! P\, -\,  \dfrac{|P|^2}{n } I\right). 
\end{align}
In view of \eqref{3.6}, formula \eqref{3.1} has been established.           \qed

\begin{lemma} \label{l1,1} Let $K$ be given by \eqref{1.3}. Then, $K\in C^2(S^+)$ and its 2nd derivative is given by
\beq \label{3.1A}
K_{PP}(P)\ =\ 2 I \ot  \frac{\big(P^\top \! P \big)^{-1} S\big(P^\top \! P \big)}{\det\big(P^\top \! P \big)^{1/n}}\ + \ 2 P \ot P :  \frac{\big(P^\top \! P \big)^{-1} E}{\det\big(P^\top \! P \big)^{1/n}}\ + \ O(P)
\eeq
which in index form can be written as
\begin{align} \label{3.2A}
K_{P_{\al i}P_{\be j}} (P)\ =\ & 2\de_{\al \be}  \Bigg( \frac{\big(P^\top \! P \big)_{ik}^{-1} \big( P_{\ga k } P_{\ga j}-\frac{1}{n}|P|^2\de_{kj} \big) }{\det\big(P^\top \! P \big)^{1/n}} \Bigg) \nonumber\\
& +\ 2P_{\al m}P_{\be l }   \Bigg( \frac{ \big(P^\top \! P \big)_{ik}^{-1}  E_{kjlm}  }{\det\big(P^\top \! P \big)^{1/n}}\Bigg) \ + \ O_{\al i \be j}(P).
\end{align}
Here $O_{\al i \be j}(P)$ is a tensor of the form $K_{P_{\al m}}(P)A_{m\be i j}(P)$ and is annihilated by $[K_P(P)]_{\ga \al}^\bot$, that is $[K_P(P)]^\bot O(P)=0$.\ $E$ is a constant 4th order tensor whose components $E_{kjlm}$ are given by
\beq  \label{3.3A}
E_{kjlm}\ := \ \de_{ml}\de_{jk} + \de_{mj}\de_{kl} -\frac{2}{n}\de_{mk}\de_{jl} .
\eeq
\end{lemma}

The explicit form of the tensor $O_{\al i \be j}(P)$ is a complicated formula which follows by the proof of Lemma \ref{l1,1}, but we do not need this formula because is ``killed" by $[K_P(P)]^\bot$ and doe not appear in \eqref{1.4}.

\BPL \ref{l1,1}. We begin by calculating the derivative $\big(\big(P^\top \! P \big)_{mi}^{-1}\big)_{P_{\be j}}$. We have
\beq \label{3.4A}
\big(P^\top \! P \big)_{mi}^{-1} \big(P^\top \! P \big)_{ik}\ = \ \de_{mk}
\eeq
which gives
\begin{align} \label{3.5A}
\big(\big(P^\top \! P \big)_{mi}^{-1} \big)_{P_{\be j}} \big(P^\top \! P \big)_{ik} \ &= \ - \big(P^\top \! P \big)_{mi}^{-1} (P_{\ga i} P_{\ga k})_{P_{\be j}} \nonumber\\
&=  \ - \big(P^\top \! P \big)_{mi}^{-1}[ \de_{\be \ga} \de_{ij}P_{\ga k} + P_{\ga i} \de_{\be \ga} \de_{kj} ]\\
&=  \ - \big(P^\top \! P \big)_{ml}^{-1} [ P_{\be k}\de_{lj} + P_{\be l} \de_{kj} ]. \nonumber
\end{align}
Hence, we have
\beq \label{3.6A}
\big(\big(P^\top \! P \big)_{mi}^{-1} \big)_{P_{\be j}}  \ =  \ - \big(P^\top \! P \big)_{ml}^{-1} [ P_{\be k}\de_{lj} + P_{\be l} \de_{kj} ] \big(P^\top \! P \big)_{ki}^{-1}. 
\eeq
Now we differentiate \eqref{3.1a}:
\begin{align} \label{3.7A}
K_{P_{\al i}P_{\be j}} (P)\ =\ \, & 2\de_{\al \be}\de_{mj}  \Bigg( \frac{ \de_{mi}-\frac{1}{n}|P|^2\big(P^\top\!P\big)^{-1}_{mi}}{\det\big(P^\top \! P \big)^{1/n}} \Bigg) \  - \ 2P_{\al m}  \Bigg( \frac{\big(|P|^2\big(P^\top\!P\big)^{-1}_{mi}\big)_{P_{\be j}}}{n\det\big(P^\top \! P \big)^{1/n}} \Bigg) 
\nonumber\\
& - \ \Bigg[  2P_{\al m}  \Bigg( \frac{ \de_{mi} -\frac{1}{n}|P|^2\big(P^\top\!P\big)^{-1}_{mi}}{\det\big(P^\top \! P \big)^{1/n}} \Bigg) \Bigg]
\frac{ \big(\det\big(P^\top \! P \big)^{1/n}\big)_{P_{\be j}} }{\det\big(P^\top \! P \big)^{1/n}}.
\end{align}
In view of \eqref{3.1a}, the last summand in \eqref{3.7A} is annihilated by the projection $[K_P(P)]_{\ga \al}^\bot$. We rewrite \eqref{3.7A} as
\begin{align} \label{3.8A}
K_{P_{\al i}P_{\be j}} (P)\ =\ \, &2\de_{\al \be}  \Bigg( \frac{ \de_{ij}-\frac{1}{n}|P|^2\big(P^\top\!P\big)^{-1}_{ij}}{\det\big(P^\top \! P \big)^{1/n}} \Bigg) \nonumber\\
&  - \ 2P_{\al m}  \Bigg( \frac{ \big(|P|^2\big(P^\top\!P\big)^{-1}_{mi}\big)_{P_{\be j}}}{n\det\big(P^\top \! P \big)^{1/n}} \Bigg) \ +\ O_{\al i \be j}(P).
\end{align}
By using \eqref{3.6A} in \eqref{3.8A}, we have
\begin{align} \label{3.9A}
K_{P_{\al i}P_{\be j}} (P)\ =& \ \, 2\de_{\al \be}  \Bigg( \frac{ \de_{ij}-\frac{1}{n}|P|^2\big(P^\top\!P\big)^{-1}_{ij}}{\det\big(P^\top \! P \big)^{1/n}} \Bigg) \ + S_{\al i \be j}(P)\ +\ O_{\al i \be j}(P),
\end{align}
where we have set
\beq \label{3.10A}
S_{\al i \be j}(P)\ := \ \frac{2}{n} P_{\al m} \frac{ 2P_{\be j} \big(P^\top\!P\big)^{-1}_{mi} -|P|^2 \big(P^\top \! P \big)_{ml}^{-1} [ P_{\be k}\de_{lj} + P_{\be l} \de_{kj} ] \big(P^\top \! P \big)_{ki}^{-1}   }{ \det\big(P^\top \! P \big)^{1/n} } .
\eeq
Equation \eqref{3.10A} gives
\begin{align} \label{3.10B}
S_{\al i \be j} (P)\ = \ & -\ \frac{4}{n} P_{\al m}P_{\be j} \frac{ \big(P^\top\!P\big)^{-1}_{mi} }{\det\big(P^\top \! P \big)^{1/n} } \nonumber\\
&+\ 2 P_{\al m} \left( \frac{\frac{1}{n}|P|^2\big(P^\top\!P\big)^{-1}_{mj}}{\det\big(P^\top \! P \big)^{1/n} }\right) \big(P^\top\!P\big)^{-1}_{ki}P_{\be k}\\
&+\ 2 P_{\al m} \left( \frac{\frac{1}{n}|P|^2\big(P^\top\!P\big)^{-1}_{mk}}{\det\big(P^\top \! P \big)^{1/n} }\right) \big(P^\top\!P\big)^{-1}_{ij}P_{\be k}. \nonumber
\end{align}
We rewrite \eqref{3.10B} as
\begin{align} \label{3.11A}
S_{\al i \be j} (P)\ = \ & -\ \frac{4}{n} P_{\al m}P_{\be j} \frac{ \big(P^\top\!P\big)^{-1}_{mi} }{\det\big(P^\top \! P \big)^{1/n} } \nonumber\\
&+\ 2 P_{\al m} \left( \frac{-\de_{mj} +\frac{1}{n}|P|^2\big(P^\top\!P\big)^{-1}_{mj}}{\det\big(P^\top \! P \big)^{1/n} } \ + \ \frac{\de_{mj} }{\det\big(P^\top \! P \big)^{1/n} } 
\right) \big(P^\top\!P\big)^{-1}_{ki}P_{\be k}\\
&+\ 2 P_{\al m} \left( \frac{-\de_{mk}+\frac{1}{n}|P|^2\big(P^\top\!P\big)^{-1}_{mk}}{\det\big(P^\top \! P \big)^{1/n} }\ +\ \frac{\de_{mk}}{\det\big(P^\top \! P \big)^{1/n} }\right) \big(P^\top\!P\big)^{-1}_{ij}P_{\be k} \nonumber
\end{align}
and observe that in view of \eqref{3.1a}, $[K_P(Du)]_{\ga \al}^\bot$ annihilates the first summands in the brackets of \eqref{3.11A} and $S_{\al i \be j} (P)$ simplifies to
\begin{align} \label{3.12A}
S_{\al i \be j} (P)\ =& \  \, 2\frac{ P_{\al k}P_{\be k}  \big(P^\top\!P\big)^{-1}_{ij} 
+ P_{\al j}P_{\be k}  \big(P^\top\!P\big)^{-1}_{ki} - \frac{2}{n} P_{\al m}P_{\be j}\big(P^\top\!P\big)^{-1}_{mi} 
}{\det\big(P^\top \! P \big)^{1/n} } \nonumber\\
& +\ O_{\al i \be j} (P),
\end{align}
for some tensor $O_{\al i \be j} (P)$ annihilated by $[K_P(Du)]_{\ga \al}^\bot$. We rewrite \eqref{3.12A} as
\begin{align} \label{3.13A}
S_{\al i \be j} (P)\, = \,  2 P_{\al m}P_{\be l}  \big(P^\top\!P\big)^{-1}_{ki}\left(\frac{ \de_{ml}\de_{jk} + \de_{mj}\de_{kl} -\frac{2}{n}\de_{mk}\de_{jl}  }{\det\big(P^\top \! P \big)^{1/n} } \right) + O _{\al i \be j} (P).
\end{align}
In view of \eqref{3.13A},  \eqref{3.10A},  \eqref{3.9A} and \eqref{3.3A},  equation  \eqref{3.2A} follows.
\qed

\ms
In view of Lemma \ref{l1}, the Euler-Lagrange system describing Optimal $p$-Quasiconformal immersions $u : \Om \sub \R^n \larrow \R^N$ is
\beq \label{3.8}
Q_p u\ := \ \Div\Big(K(Du)^{p-1} K_P(Du) \Big)\ = \  0.
\eeq
In view of \eqref{3.1}, \eqref{3.8} can be written in index form as
\beq \label{3.9}
D_i \left( \left(\frac{\tr(g)}{\det(g)^{1/n}}\right)^{p-1} \! D_k u_\al \,\frac{ g^{-1}_{km} S(g)_{mi}}{\det(g)^{1/n}}\right) \ = \  0,
\eeq
where $g=Du^\top \! Du$ is the Riemannian metric and $S$ is the Ahlfors operator of \eqref{2.6}.

\section{Derivation of the PDE System Governing Optimal $\infty$-Quasiconformal Immersions.} \label{section4}

The derivation we perform is this section can be deduced by a reworking of our results in \cite{K1, K2} and application of Lemmas \ref{l1} and \ref{l1,1} proved previously, but for the reader's convenience it is best to argue at the outset. Let $u : \Om \sub \R^n \larrow \R^N$ be an immersion in $C^2(\Om)^N$. By  distributing derivatives in \eqref{3.8}, we have
\begin{align}  \label{4.1}
 (p -1)  K^{p-2} K_{P_{\al i}}(Du) K_{P_{\be j}}(Du)   D^2_{ij}u_\be \ + \ 
K^{p-1}K_{P_{\al i} P_{\be j}}(Du) D^2_{ij}u_\be \ =  \ 0.
\end{align}
For each $x\in \Om$, $K_P\big((Du)(x)\big) : \R^n \larrow \R^N$  is a linear map. We define the orthogonal  projections
\begin{align}
[K_P(Du)]^\bot\ & :=\ \textrm{Proj}_{N((K_P(Du))^\top)}, \label{4.3}\\
[K_P(Du)]^\top\ & :=\ \textrm{Proj}_{R(K_P(Du))}, \label{4.4}
\end{align}
which are the projections on nullspace of $(K_P(Du))^\top$ and range of $K_P(Du)$ respectively. We rewrite \eqref{4.1} by applying the expansion $I = [K_P(Du)]^\bot + [K_P(Du)]^\top$ of the identity of $\R^N$ and contract the derivative in the left hand side to obtain
\begin{align} \label{4.5}
K_P(Du) D\big(K(Du)\big) \ \,+& \ \frac{K}{p -1} [K_P(Du)]^\top K_{PP}(Du):D^2u \nonumber\\
                        =& \ -\frac{K}{p -1} [K_P(Du)]^\bot K_{PP}(Du):D^2u.
\end{align}
The left hand side is a vector valued in $[K_P(Du)]^\top$ and the right hand side is a vector valued in $[K_P(Du)]^\bot$.  By orthogonality, left and right hand side vanish and actually \eqref{4.5} decouples to two systems. We rescale the right hand side of \eqref{4.5} by multiplying by $p-1$ and rearrange to obtain
\begin{align} 
K_P(Du) \ot K_P(Du):D^2u \ \,+ & \ \, K[K_P(Du)]^\bot K_{PP}(Du):D^2u \nonumber\\
                        =& \ -\frac{K(Du)}{p -1} [K_P(Du)]^\top K_{PP}(Du):D^2u.
\end{align}
We rewrite as
\begin{align} \label{4.8}
\Big( K_P \ot K_P + K[K_P]^\bot K_{PP}\Big)(Du):D^2u\ = \ -\frac{K [K_P]^\top K_{PP}}{ p -1}(Du):D^2u .
\end{align}
As $p \ri \infty$, \eqref{4.8} leads to \eqref{1.4}. 

\begin{remark}
We note that we can also remove the dilation function $K$ from the normal coefficient $ [K_P]^\bot K_{PP}$ with the renormalisation because it is strictly positive: $K(Du)\geq n >0$. We do not have this option in the case of  the general system \eqref{1.4a}, because $|H(Du)|$ may vanish.  However, when $n=2\leq N$ and $H(P)=|P|^2$, in \cite{K6} we  show that non-constant $\infty$-Harmonic maps have no interior gradient zeros: either $|Du|>0$ or $|Du|\equiv 0$.
\end{remark}

The next differential identity relates our system \eqref{1.4} with the seemingly different Aronsson PDE system of Capogna-Raich in \cite{CR}. In particular, it follows that even when $n=N$ the PDE system derived in \cite{CR} is only a projection of \eqref{1.4} along $[K_P(Du)]^\top$. Hence, the PDE system in \cite{CR} seems to fail to encapsulate all the information of optimised quasiconformal maps.

\bl \label{l6} Let $u : \Om \sub \R^n \larrow \R^n$ be a local diffeomorphism in $C^1(\Om)^n$. Then, we have the identity
\beq
K_P(Du)\ = \ -\frac{2K(Du)}{n}\left( (Du)^{-1,\top} -  n\frac{Du}{|Du|^2}\right)
\eeq
where $K$ and $K_P$ are given by \eqref{1.3} and \eqref{3.1}.
\el

\BPL \ref{l6}. By observing that for any invertible $A\in \R^n \ot \R^n$ there holds $A^{-1,\top}=A^{\top,-1}$, we have 
\beq \label{5.37}
\big(Du^\top \! Du \big)^{-1}\ = \ (Du)^{-1} (Du)^{\top,-1}\ =  \ (Du)^{-1} (Du)^{-1,\top}.
\eeq
Thus, we obtain
\begin{align} \label{5.38}
(Du)^{-1,\top} -  n\frac{Du}{|Du|^2}\ = & \ -\frac{n}{|Du|^2} \left( Du -  \frac{|Du|^2}{n} (Du)^{-1,\top} 
\right) \nonumber \\
= & \ -\frac{n}{|Du|^2} \left( Du -  \frac{|Du|^2}{n} Du(Du)^{-1} (Du)^{-1,\top} \right)  \\
= & \ -\frac{n}{|Du|^2} Du \left(I -  \frac{|Du|^2}{n} (Du)^{-1}(Du)^{-1,\top} \right). \nonumber
\end{align}
Consequently, by \eqref{5.37} and \eqref{5.38}, we obtain
\begin{align}
-\frac{|Du|^2}{n}\left( (Du)^{-1,\top} -  n\frac{Du}{|Du|^2} \right)\ = & \ Du \left( I -  \frac{|Du|^2}{n}\big(Du^\top \!Du \big)^{-1} \right)  \nonumber\\
= & \ Du\, \big(Du^\top \!Du \big)^{-1}\left( Du^\top \!Du -  \frac{|Du|^2}{n}I \right).
\end{align}
Hence, by  \eqref{3.1} and \eqref{1.3} we have
\begin{align}
-\frac{2K(Du)}{n}\left( (Du)^{-1,\top} -  n\frac{Du}{|Du|^2} \right)\
= &  \  2Du\, \big(Du^\top \!Du \big)^{-1}  \left( \frac{  Du^\top \!Du -  \frac{|Du|^2}{n}I }{\det\big(Du^\top \!Du \big)^{1/n}  }\right)  \nonumber\\
=&\  K_P(Du). 
\end{align}
The desired identity follows.                \qed

\section{Variational Structure of Optimal $\infty$-Quasiconformal Immersions.} \label{section5}

We begin by introducing a minimality notion of vector-valued Calculus of Variations in $L^\infty$ for the supremal dilation functional \eqref{1.2}. Let $u : \Om \sub \R^n \larrow \R^N$ be an immersion in $C^1(\Om)^N$. In view of \eqref{3.1}, we have the identity
\beq
K_P(Du)\ = \ \Bigg(2\frac{Du \big(Du^\top\! Du \big)^{-1}}{\det\big(Du^\top\! Du\big)^{1/n}} \Bigg) S\big(Du^\top\! Du \big).
\eeq
Generally, the rank of $K_P(Du)$ may not be constant throughout $\Om$, although by assumption $\rk(Du)= \rk(Du^\top\! Du) \equiv n$, because possibly $\rk(S(Du^\top\! Du))<n$ on certain regions of $\Om$. We set
\begin{align} \label{5.1a}
\Om_k \ :=\ \inter\Big\{ \rk\big(S(Du^\top\! Du)\big) \, = \, k \Big\}\ , \ \ \ k\, = \, 0,\, 1\, , ....\, , \, n,
\end{align}
where $``\inter"$ denotes topological interior. The $n+1$ open sets $\Om_k$ are the \emph{``phases''} of the immersion $u$. Their complement in $\Om$
\beq \label{5.1b}
\S\ :=\ \Om \set \left(\cup_0^n \Om_k\right)
\eeq
is the set of \emph{``interfaces''} and is closed in $\Om$ with empty interior. We will also need the \emph{``augmented phases"}
\begin{align} \label{5.3a}
\Om^*_k \ :=\ \Big\{ \rk \big(S(Du^\top\! Du)\big) \, = \, k \Big\}\ , \ \ \ k\, = \, 0,\, 1\, , ....\, , \, n.
\end{align}
Obviously, $\{\Om_0^*,...,\Om^*_n\}$ is a partition of $\Om$ to disjoint phases and $\S$ can be written as $\S=\cup_0^n (\Om^*_k \set \Om_k)$. The extreme cases of $\Om^*_0$ and $\Om^*_n$ are particularly important. $\Om^*_0$ is the \emph{conformality set} of the immersion and is closed in $\Om$. Hence,
\beq 
\Om^*_0 \ = \ \left\{ Du^\top \! Du =\frac{|Du|^2}{n}I\right\}.
\eeq
Similarly, by Corollary \ref{c7} that follows, if $u$ solves $K_P(Du) \ot K_P(Du) :D^2u=0$, then $\Om^*_n$ is the \emph{constant dilation set} of the immersion and coincides with $\Om_n$:
\beq 
\Om^*_n \ = \ \left\{ \frac{|Du|^2}{\det(Du^\top\! Du)^{1/n}}=const.\right\}.
\eeq
If $\Om_n$ is not connected, then the constants may differ in connected cmponents.

\begin{definition}\label{def1} Let $u : \Om \sub \R^n \larrow \R^N$ be an immersion in $C^1(\Om)^N$.

\ms \noi (i) We say that $u$ has \emph{Rank-One Locally Minimal Dilation} when for all compactly contained subdomains $D$ of $\Om$, all functions $g$ over $D$ vanishing on $\p D$ and all directions $\xi$, $u$ is a minimiser on $D$ with respect to essentially scalar variations $u+f\xi$:
\beq \label{5.2}
\left.
\begin{array}{l}
D \subset \subset \Om, \\
f\in C^1_0(D), \\
\xi \in \mS^{N-1}
\end{array}
\right\} \ \ \Longrightarrow \ \
K_\infty(u,\Om)\ \leq \ K_\infty(u+f\xi,\Om).
\eeq
\[
\underset{\text{Figure 2.}}{\includegraphics[scale=0.24]{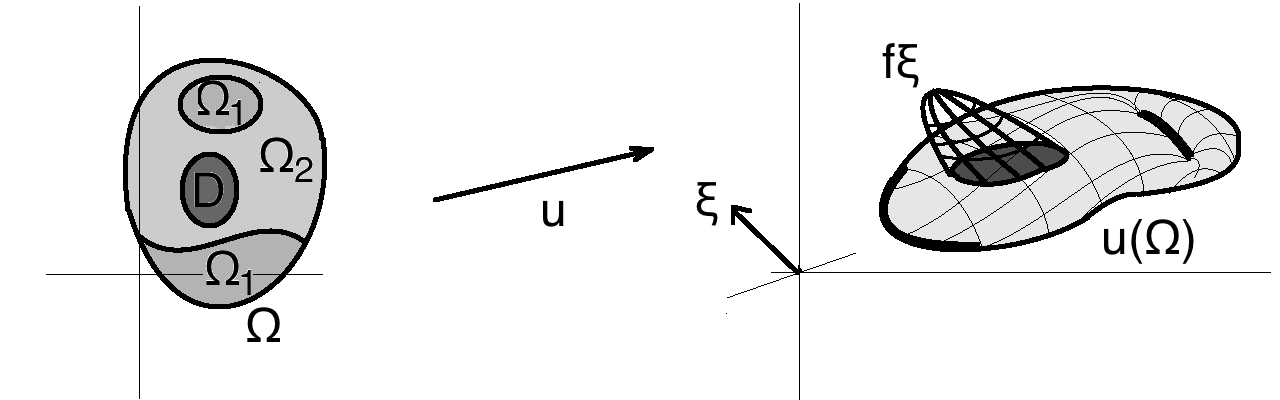}} \label{fig2} 
\]
 \noi (ii) We say that \emph{$u(\Om)$ has Minimally Distorted Area} when for all compactly contained subdomains $D$ \emph{off the interfaces}, all functions $h$ on $\bar{D}$ (not only vanishing on $\p D$) and all vector fields  $\nu$ along $u$ normal to $K_P(Du)$, $u$ is a minimiser on $D$ with respect to normal free variations $u+h\nu$:
\beq \label{5.3}
\left.
\begin{array}{l}
D \subset \subset \Om \set \S, \\
h\in C^1(\bar{D}), \\
\nu \in \Gamma([K_P(Du)]^\bot)
\end{array}
\right\} \ \ \Longrightarrow \ \
K_\infty(u,\Om)\ \leq \ K_\infty(u+h \nu,\Om).
\eeq
\[
\underset{\text{Figure 3.}}{\includegraphics[scale=0.24]{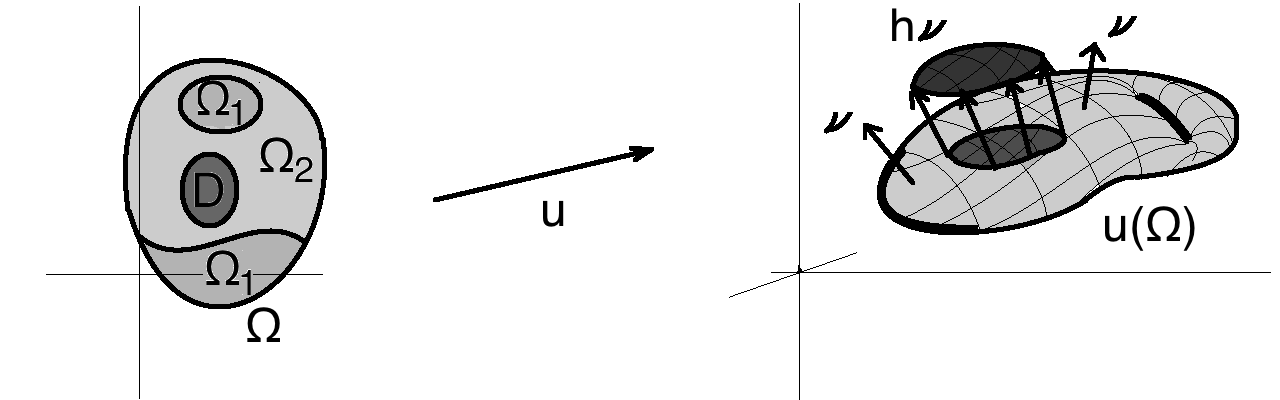}}  \label{fig3}  
\]

\noi (iii) We call $u$ \emph{Minimal $\infty$-Quasiconformal Immersion} when $u$ is has Rank-One Locally Minimal Dilation with Minimally Distorted Area of $u(\Om) \sub \R^N$.
\end{definition}

By employing the previous minimality notion, we have the next

\bt[Variational Structure of Optimal $\infty$-Quasiconformal Immersions] \label{th1}  Let $u : \Om \sub \R^n \larrow \R^N$ be an immersion in $C^2(\Om)^N$. Then, if $u$ is Minimal $\infty$-Quasiconformal, it follows that $u$ solves
\begin{align}
K_P(Du) \ot K_P(Du) :D^2u\ &= \ 0, \ \text{ on }\Om,\\
[K_P(Du)]^\bot K_{PP}(Du):D^2u\ &= \ 0,  \ \text{ on }\Om\set \S, \label{5.9a}
\end{align}  
where $\S$ is the set of interfaces of rank discontinuities of $S(Du^\top\! Du)$.
\et

We note that by the results of Section \ref{section6} that follows, in the case $n=2 \leq N$ Theorem \ref{th1} can be strengthend to the following

\begin{corollary}[2-Dimensional Optimal $\infty$-Quasiconformal Immersions] \label{cor11}  Let $u : \Om \sub \R^2 \larrow \R^N$ be an immersion in $C^2(\Om)^N$. If $u$ is Minimal $\infty$-Quasiconformal, it follows that $u$ is Optimal $\infty$-Quasiconformal.
\end{corollary}

The point in Corollary  \ref{cor11} is that \eqref{5.9a} is satisfied on $\Om$ and not only on $\Om\set \S$. Actually, \emph{when $n=2$ then the set of interfaces is empty: $\S=\emptyset$}.

\ms
The proof of Theorem \ref{th1} is split in two lemmas. 

\bl \label{l2} Let $u : \Om \sub \R^n \larrow \R^N$ be an immersion in $C^2(\Om)^N$. If $u$ has Rank-One Locally Minimal Dilation, then $u$ solves $K_P(Du) \ot K_P(Du)  :D^2u=0$ on $\Om$.
\el
The proof of Lemma \ref{l2} follows by Theorem 2.1 in \cite{K1} and relates to Lemma 2.3 in \cite{K2}, but we present a simplified more direct proof for the reader's convenience. 

\BPL \ref{l2}. Fix $x\in \Om$, $0<\e < \dist(x,\p \Om)$, $\de>0$ and $\xi \in \mS^{N-1}$. Choose $D:=\mB_\e(x)$ and $f  \in C^1_0(D)$ given by
\beq
f(z)\ :=\ \frac{1}{2}\big(\e^2-|z-x|^2\big) .
\eeq
Since $\rk(Du)=n$ on $\Om$ and $Df(z)=-(z-x)$, by restricting $\de$ sufficiently we obtain that $\rk(Du+\de \xi \ot Df)=n$ on $\mB_\e(x)$. By Taylor expansions of $K(Du)$ and $K(Du+\de \xi \ot Df)$ at $x$ we have
\begin{align} \label{5.7}
K(Du(z))\ = \ K(Du(x))\ +\ D\big(K(Du)\big)(x)^\top(z-x)\ + \ o(|z-x|),
\end{align}
as $z\ri x$, and also by using that $D^2f=-I$ and $Df(x)=0$ we have
\begin{align} \label{5.8}
K\big((Du+\de\xi \ot Df)(z)\big)\ & = \ K\big((Du+\de\xi \ot Df)(x)\big) \nonumber\\ 
& \ \ \ \ +\ D\big(K(Du+ \de\xi \ot Df)\big)(x)^\top(z-x) \ + \ o(|z-x|)\nonumber\\
             & = \ K(Du(x))\ +\ K_P(Du(x))^\top \big(D^2u(x)-\de \xi\ot I \big)(z-x) \\
     &\ \ \ \  + \ o(|z-x|)\nonumber\\
 & = \ K(Du(x))\, + \, \Big(D\big(K(Du)\big)^\top \! - \de \xi^\top K_P(Du)\Big)(x) (z-x) \nonumber\\
&\ \ \ \  + \ o(|z-x|),\nonumber
\end{align}
as $z\ri x$. By \eqref{5.7} we have the estimate
\begin{align} \label{5.9}
K_\infty \big(u,\mB_\e(x)\big)\ & \geq  \ K(Du(x)) \ +\ \max_{\{|z-x|\leq \e\}}\Big\{D\big(K(Du)\big)(x)^\top(z-x)\Big\}\  + \ o(\e) \nonumber\\
& = \ K(Du(x))\ +\ \e\big|D\big(K(Du)\big)(x)\big|\ + \ o(\e),
\end{align}
as $\e \ri 0$, and also by \eqref{5.8} we have
\begin{align} \label{5.10}
K_\infty \big(u+\de f\xi ,\mB_\e(x)\big) \ & \leq  \ K(Du(x)) \ + \max_{\{|z-x|\leq \e\}}\Big\{D\big(K(Du)\big)^\top \nonumber\\
&\ \ \ \ - \de \xi^\top K_P(Du)\big)(x)(z-x)\Big\} \  + \ o(\e) \\
& = \ K(Du(x))\ +\ \e\big|D\big(K(Du)\big) - \, \de \xi^\top K_P(Du)\big|(x)\ + \ o(\e) , \nonumber
\end{align}
as $\e \ri 0$. Then, since $u$ has rank-one locally minimal dilation, by \eqref{5.9} and \eqref{5.10} we have
\begin{align} \label{5.11}
0\ & \leq \ K_\infty \big(u+\de f\xi ,\mB_\e(x)\big) \, - \,  K_\infty \big(u,\mB_\e(x)\big) \nonumber\\
&\leq \  \e\Big(\big|D\big(K(Du)\big) -\, \de \xi^\top K_P(Du)\big|\ -\ \big|D\big(K(Du)\big)\big|\Big)(x)\ + \ o(\e),
\end{align}
as $\e \ri 0$. Suppose first $D\big(K(Du)\big)(x)=0$. Since
\begin{align} \label{5.11a}
K_P(Du) \ot K_P(Du) : D^2u  \ = \  K_P(Du) D \big(K(Du)\big) 
\end{align}
we obtain that $\big(K_P(Du) \ot K_P(Du): D^2u\big)(x) = 0$ as desired. If $D\big(K(Du)\big)(x)\neq0$, then Taylor expansion of the function 
\beq
p\ \mapsto\  \big|D\big(K(Du)\big)(x) +\, p\big| - \big|D\big(K(Du)\big)(x)\big| 
\eeq
at $p_0=0$ and evaluated at $p=-\, \de \xi^\top K_P(Du(x))$, \eqref{5.11} implies after letting $\e \ri 0$ that
\beq \label{5.12a}
0\ \leq  \ - \de \, \xi^\top K_P(Du(x))\left(\frac{D\big(K(Du)\big)}{\big|D\big(K(Du)\big)\big|} \right)(x) \ + \ o(\de).
\eeq
By letting $\de \ri 0$ in \eqref{5.12a} we obtain $\xi^\top \big(K_P(Du) \ot K_P(Du): D^2u\big)(x) \geq 0$ for any direction $\xi$. Since $\xi$ and $x$ are arbitrary we get $K_P(Du) \ot K_P(Du): D^2u = 0$ on $\Om$. The lemma follows.               \qed

\bl \label{l3} Let $u : \Om \sub \R^n \larrow \R^N$ be an immersion in $C^2(\Om)^N$ with Minimally Distorted Area of $u(\Om)$. Then, $u$ solves $[K_P(Du)]^\bot K_{PP}(Du) : D^2u=0$ on $\Om \set \S$.
\el

\BPL \ref{l3}. Fix $x\in \Om \set \S$. Then, $x$ belongs to some phase $\Om_k$ of constant rank and $\rk\big(S(Du^\top\! Du)\big)\equiv k$ thereon. We choose $0<\e<\frac{1}{2}\dist(x,\p \Om_k)$ and $0<\de<1$. By the Rank Theorem (see e.g.\ \cite{N}) and application of the Gram-Schmidt procedure to a local frame field adapted to the immersion near $u(x)$, we can construct a local frame of sections $\{\nu^1,...,\nu^{N-k}\}$ spanning $\Gamma([K_P(Du)]^\bot,\mB_{2\e}(x))$ for $\e$ small enough. Let  $\nu$ be a linear combination of these sections and choose an $h \in C^1\big(\overline{\mB_\e(x)} \big)$. Since $\rk(Du)=n$ on $\Om$, by restricting $\de$ sufficiently we obtain $\rk\big(D(u+\de h\nu)\big)=n$ on $\mB_\e(x)$. By differentiating $\nu^\top K_P(Du) = 0$ we obtain
\beq \label{5.13}
D_k\nu_\al K_{P_{\al i}}(Du)\ =  \ -\nu_\al K_{P_{\al i}P_{\be j}}(Du) D_{kj}^2u_{\be}
\eeq
and by putting $i=k$ and summing, we get
\beq \label{5.14}
D_i\nu_\al K_{P_{\al i}}(Du)\ =  \ -\nu_\al K_{P_{\al i}P_{\be j}}(Du) D_{ij}^2u_{\be}
\eeq
that is 
\beq \label{5.14a}
D\nu : K_P(Du) \ =\ -\nu^\top K_{PP}(Du) : D^2u. 
\eeq
By Taylor expansion of the dilation and usage of $\nu^\top K_P(Du)=0$, we obtain
\begin{align} \label{5.16}
K\big(D(u+\de h\nu)\big)\ &=\ K(Du)\ + \ K_P\big(Du):D(\de h\nu) \ + \ o(\de|h\nu|) \nonumber\\
             &=\ K(Du)\ + \ \de K_P\big(Du): \big(hD\nu \, +\, \nu \ot Dh \big) \ + \ o(\de) \\
&=\ K(Du)\ + \ \de \Big(hD\nu: K_P\big(Du) \, +\, \nu^\top K_P\big(Du) Dh\Big)\ + \ o(\de)  \nonumber\\
&=\ K(Du)\ + \ \de hD\nu: K_P\big(Du)\ + \ o(\de)  \nonumber
\end{align}
as $\de \ri 0$. By \eqref{5.16} and \eqref{5.14a} we have
\beq \label{5.18}
K\big(D(u+\de h\nu)\big)\ = \  \ K(Du) \ - \ 2\de h\big(\nu^\top K_{PP}(Du) : D^2u\big) \   + \ o(\de)  ,
\eeq
as $\de \ri 0$. Hence, since $u(\Om)$ has minimally distorted area, by \eqref{5.18} we have
\begin{align} \label{5.19}
K_\infty \big( u, \mB_\e(x)\big)\ &\leq \ K_\infty \big( u+\de h\nu, \mB_\e(x)\big) \nonumber\\
      &=\ \sup_{\mB_\e(x)} \Big\{ K(Du) \ - \ 2\de h \big(\nu^\top K_{PP}(Du) : D^2u\big)\   + \ o(\de) \Big\}
  \end{align}
as $\de \ri 0$, which gives 
\begin{align} \label{5.20}
K_\infty \big( u, \mB_\e(x)\big)\ &\leq \sup_{\mB_\e(x)} K(Du)\ -\ 2\de\min_{\overline{\mB_\e(x)}} \Big\{ h \big(\nu^\top K_{PP}(Du) : D^2u\big) \Big\}  \   + \ o(\de) \nonumber\\
  &= K_\infty \big( u, \mB_\e(x)\big)\ -\ 2\de\min_{\overline{\mB_\e(x)}} \Big\{ h \big(\nu^\top K_{PP}(Du) : D^2u\big) \Big\}  \   + \ o(\de) . 
           \end{align}
Hence, by passing to the limit as $\de \ri 0$, \eqref{5.20} gives
\beq \label{5.21}
\min_{\overline{\mB_\e(x)}} \Big\{  h\big(\nu^\top K_{PP}(Du) : D^2u\big) \Big\}\ \leq \ 0.
\eeq
We now choose as $h$ the constant function 
\beq
h\ :=\ \sgn \left( \nu^\top K_{PP}(Du) : D^2u \right)(x) 
\eeq
and by \eqref{5.21} as $\e \ri 0$ we get $\big|\nu^\top K_{PP}(Du) : D^2u\big|(x)=0$. Since $\nu$ is an arbitrary normal section and $x$ is an arbitrary point on $\Om \set \S$, we get $([K_P]^\bot K_{PP})(Du) : D^2u=0$ on $\Om \set \S$ and the lemma follows.                              \qed

\section{Geometric Properties of Optimal $\infty$-Quasiconformal Immersions.} \label{section6}

\subsection{Geometric Form of the PDE System.} \label{subsection6.1} In this subsection we show that system \eqref{1.1} decouples to two system one normal to to other which can be written in geometric rather coordinate-free fashion, at least within the phases of solutions whereon the coefficients of the system are continuous. 

\begin{proposition} \label{l5} Let $u : \Om \sub \R^n \larrow \R^N$ be an  immersion in $C^2(\Om)^N$. If $K$ is the dilation \eqref{1.3} and its derivatives are given by \eqref{3.1} and \eqref{3.1A}, then the Aronsson system
\beq \label{5.25}
Q_\infty u \ =\ \Big(K_P \ot K_P + [K_P]^\bot K_{PP}\Big)(Du):D^2u\ = \ 0
\eeq
is equivalent on each phase $\Om_k = \inter\{\rk(S(Du^\top\! Du)) = k\}$ to the pair of systems
\begin{align}
S(\bold{G})D\big(\tr(\bold{G})\big)\ = \ 0, \label{5.26}\\
\mathbb{B}^\bot : \big(\tr(\bold{G})\big)_{P} \ =  \ 0, \label{5.27}
\end{align}
where $\bold{G}$ is given by \eqref{1.1}, $g=Du^\top\! Du$ is the Riemannian metric on $u(\Om)$, $S$ is the Ahlfors operator and $\mB^\bot$ is the ``generalized 2nd fundamental form'', defined for every local normal section $\nu \in \Gamma([K_P(Du)]^\bot,D)$ over $D\sub \Om\set \S$ as $(\mB^\bot) _\nu:= D\nu$. Moreover, \eqref{5.26} is valid on all of $\Om$.
\end{proposition}

We observe that system \eqref{5.26} can also be written as
\beq
S(g)D\left(\frac{\tr(g)}{\det(g)^{1/n}}\right)\ =  \ 0 
\eeq
and hence depends only on the metric structure of the immersion. System \eqref{5.26} is the ``tangential system". On the other hand, \eqref{5.27} can be written also as
\beq
\mathbb{B}^\bot :\left(\frac{\tr(g)}{\det(g)^{1/n}}\right)_P \ =  \ 0 
\eeq
and depends on the exterior geometry as well, the ``shape" of $u(\Om)$. System \eqref{5.27} is the ``normal system".

\BPP \ref{l5}. By applying the orthogonal projections \eqref{4.3} and \eqref{4.4} to \eqref{5.25}, we decouple it to
\begin{align}
K_P(Du) \ot K_P(Du) :D^2u\ = \ 0, \label{5.30}\\
[K_P(Du)]^\bot K_{PP}(Du) :D^2u\ =  \ 0. \label{5.31}
\end{align}
In view of \eqref{3.1}, we rewrite \eqref{5.30} as
\beq
Dug^{-1}S(g)D\big(K(Du)\big)\ = \ 0.
\eeq
By using that $K(Du)=\tr(\bold{G})$ and that $Dug^{-1}$ has constant rank equal to $n$ and hence is left invertible, we obtain
\beq \label{5.33}
\big(Dug^{-1}\big)^{-1}Dug^{-1}S(g)D\big(\tr(\bold{G})\big)\ = \ S(g)D\big(\tr(\bold{G})\big)\ = \ 0.
\eeq
Since $g=\det(g)^{1/n}\bold{G}$, system \eqref{5.33} leads to \eqref{5.26}. To obtain \eqref{5.27}, we observe that \eqref{5.31} is equivalent to
\beq \label{5.34}
\nu^\top K_{PP}(Du) :D^2u\ = \ 0,
\eeq
for all local normal sections $\nu \in \Gamma([K_P(Du)]^\bot,D)$, $D\sub \Om \set \S$. By \eqref{5.14a}, equation \eqref{5.34} is equivalent to $- D\nu:K_P(Du)= 0$. Hence, we rewrite it as
\beq \label{5.35}
- D\nu:\big(\tr(\bold{G})\big)_P\ = \ 0.
\eeq
By definition of $\mB^\bot$, system \eqref{5.35} leads to \eqref{5.27} and the proposition follows.        \qed

\begin{remark} We will later show that the 2-dimensional case $n=2\leq N$ is prominent. In this case, interfaces of discontinuities of the coefficients disappear and $\mB^\bot$ conicides with the standard 2nd fundamental form.
\end{remark}

\begin{corollary}[Constant dilation on $\Om_n$] \label{c7} Let $u :\Om \sub \R^n \larrow \R^N$ be an immersion in $C^2(\Om)^N$ solving $K_P(Du) \ot K_P(Du):D^2u=0$. Then, on the $n$-phase $\Om_n$ given by \eqref{5.1a}, $u$ has constant dilation on each connected component of  $\Om_n$.
\end{corollary}

\BPCOR \ref{c7}. By \eqref{5.1a} and \eqref{5.33}, we have that $S(g)$ is invertible on $\Om_n$ and consequently we get $D\big(K(Du)\big)=0$ on $\Om_n$.
\qed

\subsection{A Geometric Property of Interfaces of Solutions.} \label{subsection6.2} We begin with a differential identity valid \emph{on the interfaces} of discontinuity, under a local regularity assumption on the interface. We assume only $C^1$ regularity, but we allow for possibly complicated topology and self-intersections.

\begin{proposition}[Covariant Derivatives on Interfaces] \label{pr3} Let $u : \Om \sub \R^n \larrow \R^N$ be an  immersion in $C^2(\Om)^N$. Suppose the set of interfaces $\S$ inside $\Om$ given by \eqref{5.1b} contains a $C^1$ immersed submanifold $M$ and let $\nabla^M$ be its Riemannian gradient. Then, we have the identity
\begin{align} \label{5.41}
\nabla^M\big([K_P(Du)]^\bot\big): K_P(Du)\ =& \  -\big([K_P]^\bot  K_{PP}\big)(Du):D^2u\nonumber\\
&\ +\  \big([K_P]^\bot K_{PP}\big)(Du) : \nabla^{M^\bot}\! Du ,
\end{align}
valid on $M\sub \S$, where $\nabla^{M^\bot}\! $ is the orthogonal complement of $\nabla^M$ in $\R^n$.
\end{proposition}

\[
\underset{\text{Figure 4.}}{\includegraphics[scale=0.24]{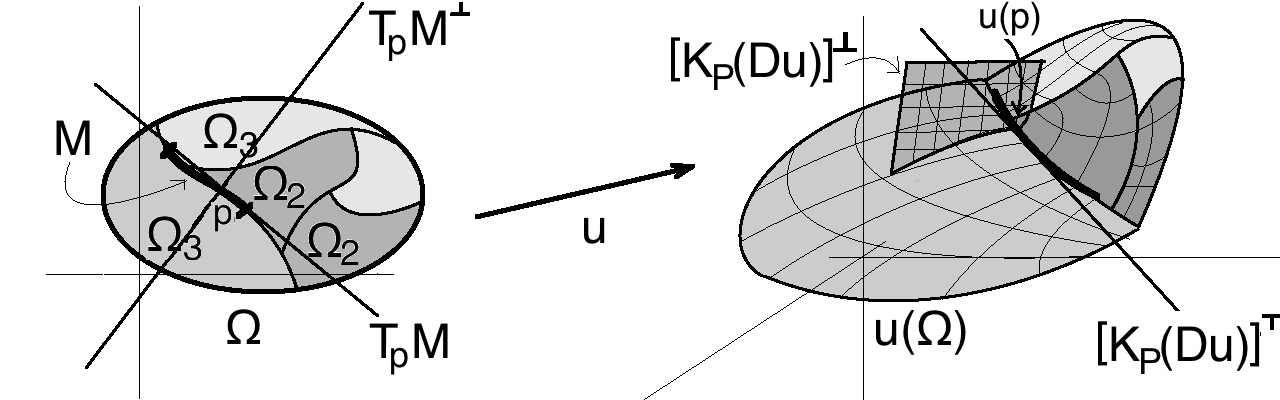}} \label{fig4} 
\]

\begin{remark} The point in \eqref{5.41} is that $[K_P(Du)]^\bot$ has covariantly differentiable contraction with $K_P(Du)$ along (part of the interface) $M$, \emph{without having assumed that $S(Du^\top\! Du)$ has constant rank on $M$} and hence without having assumed that $[K_P(Du)]^\bot$ is differentiable on $M\sub \Om$.  
\end{remark}

\BPP \ref{pr3}. By assuming as we can that $M$ is immersed by the inclusion into $\Om$, we fix a point $p \in M \sub \Om$ and consider coordinates near $p$ adapted to the immersion. Let $\{\nabla^M_1,...,\nabla^M_n\}$ denote the $n$ components of $\nabla^M$ with respect to the standard coordinates of $\R^n$. By differentiating covariantly near $p$ the identity 
\beq
[K_P(Du)]^\bot_{\al \be} K_{P_{\be j}}(Du)\ =\ 0
\eeq
we obtain
\begin{align} \label{5.43}
\nabla^M_i\big( [K_P(Du)]_{\al \be}^\bot\big)K_{P_{\be j}}(Du)\ &=- \ [K_P(Du)]_{\al \be}^\bot \nabla^M_i \big(K_{P_{\be j}}(Du)\big) \nonumber\\
&=- \ [K_P(Du)]_{\al \be}^\bot  K_{P_{\be j}P_{\ga k}}(Du)\nabla^M_i D_ku_\ga.
\end{align}
By applying the expansion $\nabla^M = D - \nabla^{M^\bot}$, putting $i=j$ and summing, \eqref{5.43} implies \eqref{5.41} and the proposition follows.                \qed

The previous identity readily implies the next

\begin{corollary} \label{c6} In the setting of Proposition \ref{pr3} above, if $u$ solves the system $([K_P]^\bot K_{PP})(Du) :D^2u=0$, then we have
\beq \label{5.45}
\nabla^M\big([K_P(Du)]^\bot\big): K_P(Du) \ =\ \big( [K_P]^\bot K_{PP}\big)(Du) :\nabla^{M^\bot}\! Du.
\eeq
In particular, the vector field 
\beq
\nabla^M\big([K_P(Du)]^\bot\big):K_P(Du)\ :\ M \larrow \R^N 
\eeq
is ``normal"to $u(M)$, namely, it is valued in $[K_P(Du)]^\bot$:
\beq \label{6.18}
[K_P(Du)]^\top \Big(\nabla^M\big([K_P(Du)]^\bot\big):K_P(Du)\Big) \ = \ 0.
\eeq
\end{corollary}

\BPCOR \ref{c6}. Since the immersion $u$ solves $[K_P(Du)]^\bot K_{PP}(Du):D^2u=0$, \eqref{5.41} gives \eqref{5.45}. By applying the projection $[K_P(Du)]^\top$ to the latter, \eqref{6.18} follows. Hence, the vector field $\nabla^M\big([K_P(Du)]^\bot\big):K_P(Du)$ equals its projection on $[K_P(Du)]^\bot$ and the corollary follows.         \qed

\ms

\section{Sufficiency of $K_P(Du) \ot K_P(Du) : D^2u =0$ for Rank-One Locally Minimal Dilation When $n=2\leq N$.} \label{section7}

In this section we show that in the case of 2-dimensional immersions when $n=2\leq N$, the tangential  system $K_P(Du) \ot K_P(Du) :D^2u=0$ is sufficient for the minimality notion of Rank-One Locally Minimal Dilation. This follows as a corollary of the fact that when $n=2$, solutions of this system necessarily have constant dilation. In particular, the rank of $S(Du^\top \!Du)$ is constant throughout the domain and interfaces of discontinuity on the coefficents of the normal system $([K_P]^\bot K_{PP})(Du):D^2u=0$ disappear.

As a corollary, we show that when $n=N=2$, the conjecture of Capogna-Raich in \cite{CR} on the sufficiency of system $(K_P\ot K_P)(Du) :D^2u=0$ for their stronger local minimality notion is false. This follows by Example \ref{ex1} below in which we construct a diffeomorphism with constant dilation on a domain of the plane which has the same boundary values with the identity.

\begin{lemma}[Constant dilation] \label{p7} Let $u : \Om \sub \R^n \larrow \R^N$ be an immersion in $C^2(\Om)^N$ which solves $K_P(Du) \ot K_P(Du) :D^2u=0$ on $\Om$. Suppose $\Om$ is connected and let $\Om^*_0,..., \Om^*_n$ be the augmented $n+1$ phases of the immersion given by \eqref{5.3a}. Then:

\noi (i) $S(Du^\top\!Du)$ has nowhere rank equal to one:
\beq
\Om^*_1\, =\ \emptyset.
\eeq
\noi (ii) 
If moreover $n=2$, then $\Om^*_0 \in \{\emptyset, \Om\}$. That is, $\Om^*_0$ is either empty or equals the whole $\Om$. Hence, $u$ has constant dilation everywhere on $\Om$:
\beq
K(Du)\  \equiv \ k \ \geq \ 2.
\eeq
If it happens that $\Om^*_0\neq \emptyset$, then $k=2$ and in this case $u$ is conformal on $\Om$.
\end{lemma}

\BPL \ref{p7}. $(i)$ On $\Om^*_1$ we have $\rk(S(Du^\top \!Du))=1$ and also $S(Du^\top \!Du)=S(Du^\top \!Du)^\top$. Since $S(Du^\top \!Du)$ is a rank-one symmetric matrix, there exist $\la : \Om^*_1 \larrow \R$ and $a : \Om^*_1 \larrow \R^n$ such that $\la>0$, $|a|=1$ and $S(Du^\top \!Du)= \la\, a \ot a$. Hence, we obtain
\begin{align}
\la\ = \ \la \, |a|^2 \ = \ \tr(\la \, a\ot a)\ =\ \tr \big(S(Du^\top \!Du)\big)\ = \  0.
\end{align}
Consequently, $\Om^*_1 = \emptyset$. 

$(ii)$ When $n=2$, by $(i)$ we have that $\Om =\Om^*_0 \cup \Om^*_2$. On $\Om^*_0$ the immersion $u$ is conformal. By Corollary \ref{c7}, on $\Om^*_2$ $u$ has constant dilation. Hence, $u$ has constant dilation on each connected component of $\Om^*_0 \cup \Om^*_2=\Om$. This means that $K(Du)$ is piecewise constant on $\Om$. By assumption, $\Om$ is connected and also $K(Du) \in C^0(\Om)$. As a result, necessarily either $\Om^*_0 =\emptyset$ or $\Om^*_0=\Om$. If $\Om^*_0 \neq \emptyset$, then $u$ is conformal on $\Om$. The lemma follows.         \qed

\ms

\begin{proposition}[Equivalences in the 2-Dimensional case] \label{c8} Let $u : \Om \sub \R^2 \larrow \R^N$ be an immersion in $C^2(\Om)^N$. Then, the following are equivalent:

\ms

\noi (i) $u$ has Rank-One Locally Minimal Dilation on $\Om$.

\ms

\noi (ii) $u$ solves $K_P(Du) \ot K_P(Du) :D^2u=0$ on $\Om$.

\ms

\noi (iii) $u$ has constant dilation on connected components of $\Om.$
\end{proposition}

\BPP \ref{c8}. The implications $(i)\Rightarrow (ii)$ and $(ii)\Rightarrow (iii)$ have already been estabished, so it suffices to prove  $(iii)\Rightarrow (i)$. For, suppose $u$ has constant dilation on connected components of $\Om$. Fix $D\subset \subset \Om$, $f\in C^1_0(D)$ and $\xi \in \mS^{N-1}$. We may assume $D$ is connected and that $\rk(Du+\xi \ot Df)=n$ on $D$. Then, since $f|_{\p D}\equiv 0$, there exists an interior critical point $\bar{x}\in D$ of $f$. By using that $Df(\bar{x})=0$, we estimate
\begin{align}
K_\infty (u+f\xi, D)\ &= \ \sup_D K\big(Du \, +\, \xi \ot Df\big) \nonumber\\
            & \geq \ K\big(Du(\bar{x})\, +\, \xi\ot Df(\bar{x})\big) \nonumber\\
            &=\ K(Du(\bar{x})) \\
&=\ \sup_D K(Du) \nonumber\\
&=\ K_\infty(u,D). \nonumber
\end{align}
Hence, $u$ has rank-one locally minimal dilation and the proposition follows.
\qed

\ms

Directly from Proposition \ref{c8} we obtain the following

\begin{corollary}[Absence of Interfaces in the 2-Dimensional case] \label{c9} Let $u : \Om \sub \R^2 \larrow \R^N$ be an immersion in $C^2(\Om)^N$ which  solves $Q_\infty u=0$ on the connected set $\Om$. Then the rank of $S(Du^\top\!Du)$ is constant on $\Om$, and equals either 0 or 2. If $rk\big(S(Du^\top\!Du)\big)=0$ then $u$ satisfies 
\begin{align}
&K(Du)\ \equiv \ 2, \\
&K_{PP}(Du) :D^2u\ = \  0.
\end{align}
The condition $K(Du) \equiv 2$ is equivalent to Conformality: $Du^\top\! Du=\frac{1}{n}|Du|^2I $. If $rk\big(S(Du^\top\!Du)\big)=2$, then $u$ satisfies
\begin{align}
&K(Du)\ \equiv \ const.\ > \ 2, \\
&[Du]^\bot K_{PP}(Du) :D^2u\ = \  0.
\end{align}
\end{corollary}

\begin{remark} Since the dilation \eqref{1.3} fails to be convex, it seems that sufficiency of the normal system $[K_P(Du)]^\bot K_{PP}(Du) :D^2u=0$ for minimally distorted area does not hold. In particular, the respective convexity arguments used in the case of the $\infty$-Laplacian in \cite{K2} fail.
\end{remark}

The following example certifies that the variational notion of rank-one locally minimal dilation is genuinely weaker than the respective notion of ``locally minimal dilation" used in \cite{CR}, where general vector-valued variations with the same boundary values are considered.

\begin{example}[Rank-One Locally Minimal Dilation is Strictly Weaker Notion] \label{ex1} (cf. \cite{CR}, Cor 1.6(2))
Let $\Om:= \mathbb{D}^2\set\{0\}\sub \R^2$ be the punctured unit disc on the plane. Fix $\ga>-1$ and consider the maps $u,u^\ga : \Om \larrow \Om$ where $u(x):=x$ and $u^\ga(x):=|x|^\ga x$. Then, $u=u^\ga$ on $\p \Om = \mS^1 \cup \{0\}$ and $u$ is conformal on $\Om$ while $u^\ga$ is quasiconformal but has constant strictly greater dilation: 
\beq
K(Du)\ \equiv \ 2 \ < \ 2\, +\, \frac{\ga^2}{\ga +1}\ \equiv \ K(Du^\ga).
\eeq
For completeness, we provide some details of our calculations. We readily have 
\beq
Du^\ga(x)\ =\ |x|^\ga\Big(I+\ga \frac{x}{|x|} \ot \frac{x}{|x|} \Big)
\eeq
and by setting $\frac{x}{|x|}=(a,b)^\top$ we obtain
\beq
Du^\ga(x)\ = \ |x|^\ga
\left[
\begin{array}{cc}
1+\ga a^2 & \ga a b\\
\ga ba & 1+\ga b^2 
\end{array}
\right]  .
\eeq
By using that $a^2 + b^2=1$, we have
\begin{align}
K(Du^\ga)\ &= \ \frac{|Du^\ga|^2}{(\det(Du^\ga)\det(Du^\ga))^{1/2}} \nonumber\\
       &=\ \frac{|x|^{2\ga} \big[(1\, +\, \ga a^2)^2 \, +\,  (1+\ga b^2)^2 \,+ \,  2(\ga ab)^2 \big]}{|x|^{2\ga} \big[ (1\, +\, \ga a^2)(1+\ga b^2)-(\ga ab)^2  \big]}\\
&=  \ 2\, +\, \frac{\ga^2}{\ga +1} . \nonumber
\end{align}
As a conclusion, in view of Corollary \ref{c8}, $u^\ga$ has rank-one minimal dilation over $\Om$, but does not have minimal dilation over $\Om$ since it has the same boundary values on $\p \Om$ with a conformal map. If moreover $\ga >0$, then both $u,u^\ga$ are in $C^1(\overline{\Om})^2$. 
\end{example}

\subsection{On the sufficiency of $K_P(Du)\ot K_P(Du) :D^2u=0$ for rank-one locally minimal dilation in the case of dimensions $3\leq n\leq N$.}

In this subsection we loosely discuss the much more complicated case of dimensions $n\geq 3$. In this case results are less sharp since Lemma \ref{p7} generally fails when $n>2$. 

To begin with, let $u : \Om \sub \R^3 \larrow \R^N$ be an immersion in $C^2(\Om)^N$. Obviously, we have $\rk(Du)=3\leq N$. By Lemma \ref{l1} and Proposition \ref{l5}, we may rewrite system $K_P(Du)\ot K_P(Du) :D^2u=0$ as
\beq \label{6.15}
g^{-1}S(g)D\big(K(Du)\big)\ = \  0,
\eeq
where $g= Du^\top\! Du$. We recall that in the case of $n=2$, Lemma \ref{p7} asserts that $S(g)$ either has two nonzero opposite eigenvalues (and hence has a saddle structure), or it vanishes. In the two-dimensional case  this covers all possible values of rank and it follows that the dilation is constant throughout connected domains. 

When $n=3$, Lemma \ref{p7} still works with the same proof, but now asserts only that 
\ms

\noi(i) there is no one-dimensional phase $\Om^*_1$, and

\ms

\noi (ii) $\Om = \Om^*_0 \cup \Om^*_2 \cup \Om^*_3$ with $K(Du)$ constant on connected components of the set $\Om^*_0 \cup \Om^*_3$. 

\ms

\noi When $n=3$ no information is provided for the two-dimensional phase $\Om^*_2$. Let us analyse more closely what happens in this case when $\Om^*_2\neq \emptyset$ and nontrivial interfaces of discontinuities may appear, where $\Om^*_2 = \{\rk(S(g))=2\}$. Let $0< \la_1 \leq \la_2 \leq \la_3 $ be the eigenvalue functions on $\Om$ of the Riemannian metric $g$. Then, the spectrum of $S(g)$ is
\begin{align} \label{6.16}
\si\big(S(g)\big)\ &=\ \si(g)\ -  \ \frac{\tr(g)}{3}  \nonumber\\
&= \ \left\{\la_1- \frac{\la_1 +\la_2 +\la_3}{3},\, \la_2- \frac{\la_1 +\la_2 +\la_3}{3},\, \la_3- \frac{\la_1 +\la_2 +\la_3}{3} \right\} \\
&= \ \left\{  \frac{2\la_1 -\la_2 -\la_3}{3},\,  \frac{2 \la_2 - \la_3 - \la_1}{3},\,  \frac{2\la_3 -\la_2 -\la_1}{3} \right\}. \nonumber
\end{align}
We distinguish the following cases:

(a) $0<\la_1=\la_2=\la_3=:\la$. Then, by \eqref{6.16} we have that $S(g)=0$.

(b) $0<\la_1=\la_2=:\la  < \la_3$. Then, by \eqref{6.16} we have that 
\beq
\si\big(S(g)\big)\ =\ \{-\mu,-\mu,2\mu\}
\eeq
where $\mu := \frac{\la_3 - \la}{3}>0$. By the Spectral Theorem, there is an orthonormal frame $\{a_1,a_2,a_3\}$ of $\R^3$ such that
\beq
S(g)\ =  \ -\mu \big(a_1 \ot a_1 \, +\, a_2\ot a_2\big) \ +\ 2\mu \, a_3 \ot a_3
\eeq
and $S(g)$ has rank three. 

(c) $0<\la_1<\la_2 = \la_3$. Again as before  $S(g)$ has rank three.

(d) $0<\la_1<\la_2 < \la_3$. This is the only case where rank equal to two may appear. Since $\la_2 +\la_3 >2\la_1$ and $\la_1 +\la_2<2\la_3$, we get
\beq
\mu_1\, :=\,  \frac{2\la_1 -\la_2 -\la_3}{3}\, <\, 0\ , \ \ \ \mu_3\, :=\,  \frac{2\la_3 -\la_2 -\la_1}{3} \, >\, 0
\eeq
but it may happen that
\beq
\mu_2\, :=\, \frac{2 \la_2 - \la_3 - \la_1}{3}
\eeq
vanishes, like for example in the extremal quasiconformal map $u : \R^3  \larrow \R^3$ given by $u(x,y,z):=(e^x,\sqrt{2}y e^x, \sqrt{3}z e^x)^\top$. We have
\beq
Du^\top\! Du\, (x,y,z) \ =  \ e^{2x}\left[
\begin{array}{ccc}1 & 0 & 0\\
0 & 2 & 0 \\
0 & 0 & 3
\end{array}
\right]
\eeq
and hence we get $(\la_1,\la_2,\la_3)=(e^{2x},2e^{2x},3e^{2x})$, which implies $\mu_2=0$. Generally, the set of interfaces of a three-dimensional optimal quasiconformal map is given by
\beq
\S\ =\ \p \{ \mu_2 \, =\, 0\}
\eeq
and the two-dimensional phase of $u$ is  given by
\beq
\Om_2\ =\ \inter\{\mu_2 = 0\} .
\eeq
Since $S(g)$ is traceless, the condition $\tr(S(g))=0$ implies $-\mu_1=\mu_3=:\mu >0$ and hence $\si\big(S(g)\big)=\{-\mu,0,\mu\}$. By the Spectral Theorem, there exists an orthonormal frame $\{a,b,c\}$ of $\R^3$ such that
\beq
S(g)\ = \ - \mu\, \big(a \ot a \ -\  c\ot c\big).
\eeq
By \eqref{6.15}, we have that $D \big(K(Du)\big)$ is perpendicular to $\{a,c\}$ and hence
\beq
D \big(K(Du)\big)\ = \ b\ot b \, D \big(K(Du)\big)
\eeq
which implies that the dilation of $u$ varies only in the direction of $b$. Consequently, $K(Du)$ depends only on $b$ through a certain function $k$:
\beq
K\big(Du(x)\big)\ = \ k\big(b(x)\big).
\eeq
Unlike the case $n=2$, when $n=3$ we do \emph{not} obtain that the dilation of three-dimensional optimal quasiconformal immersions is constant, at least not by the previous reasoning. 

However, by Theorem \ref{th1} in all dimensions $2\leq n \leq N$ rank-one locally minimal dilation implies solvability of $K_P(Du)\ot K_P(Du) :D^2u=0$ and by the higher-dimensional extension of Example \ref{ex1},  rank-one locally minimal dilation is genuinely weaker than locally minimal dilation. Although it seems reasonable that $K_P(Du)\ot K_P(Du):D^2u=0$ is sufficient for rank-one locally minimal dilation, we  can not definitely conclude for the validity of the conjecture of Capogna-Raich in \cite{CR} for $n\geq 3$.

\ms

\noi \textbf{Acknowledgement.} I am indebted to L.\ Capogna, J.\ Manfredi and Y.\ Yu for their interest in the author's work, their encouragenment and their constructive suggestions.

\bibliographystyle{amsplain}

\end{document}